\documentclass[12pt,a4paper]{amsart}
\usepackage[latin1]{inputenc}
\usepackage{latexsym}
\usepackage{color,graphicx,shortvrb}
\usepackage{amsmath, amssymb}
\usepackage{amsfonts}
\usepackage[colorlinks, bookmarks=true]{hyperref}

\newtheorem{theorem}{Theorem}[section]

\theoremstyle{definition}
\newtheorem{definition}[theorem]{Definition}
\newtheorem{remark}[theorem]{Remark}

\newcommand{\field}[1] {\mathbb{#1}}
\newcommand{\R}{\field{R}}

\newcommand{\N}{\mathbb{N}}
\newcommand{\C}{\mathbb{C}}

\newcommand{\D}{\mathbb{D}}

\title{A lethargy result for real analytic functions}
\author{J. M. Almira}

\begin{document}

 \baselineskip=16pt

\numberwithin{equation}{section}

\maketitle \markboth{Lethargy result for real analytic functions}{J. M. Almira}

\begin{abstract}
We prove that, if $(C[a,b],\{A_n\})$ is an approximation scheme and $(A_n)$ satisfies de La Vallée-Poussin Theorem, there are instances of continuous functions on $[a,b]$,
real analytic on $(a,b]$, which are ``poorly approximable'' by the elements of  $\{A_n\}$. This i\-llus\-tra\-tes the thesis
that the smoothness conditions guaranteeing that a function is ``well approximable''
must be ``global''. The failure of smoothness at endpoints may result in an
arbitrarily slow rate of approximation. 
\end{abstract}

\section{Motivation} After the work by Bernstein  \cite{bernsteininverso} and  Jackson \cite{jackson} in Approximation Theory it is accepted that there is a strong relationship between the smoothness properties of a function and the rate of decay to zero of the sequences of errors of best uniform polynomial approximation $\{E(f,\Pi_n)\}$. Indeed, this connection between smoothness and the order of convergence to zero of a sequence of best approximation errors holds for many other approximation schemes such as rational approximation, spline approximation, etc. There are many papers and monographs devoted to show this claim in its many distinct cases (see, for example, \cite[Chapter 7]{devore} and \cite{pinkus}, just to mention two beautiful references where this philosophy is neatly shown) In particular, if the approximation scheme $(X,\{A_n\})$ satisfies both Jackson's inequality $E(f,A_n)\leq Cn^{-r}|f|_Y$, $(f\in Y)$, and Bernstein's inequality $|a_n|_Y\leq Cn^r\|a_n\|_X$, $(a_n\in A_n)$, for $n=0,1,\cdots $, with respect to a  quasi-semi-normed proper subspace $Y$ of $X$, the approximation space $A_q^{\alpha}(X,\{A_n\})=\{x\in X:\{(n+1)^{\alpha-\frac{1}{q}}E(x,A_n)\}\in \ell_q\}$ coincides with the interpolation space $(X,Y)_{\alpha/r,q}$ for $0<\alpha<r$ and $0<q\leq\infty$, which is usually viewed as a smoothness space (see \cite[Chapter 7, Theorem 9.1]{devore} for a proof of this result). Another well known result which is usually understood as a proof of the strong relationship that exists between smoothness of functions and the order of decay to zero of best approximation errors is Bernstein's characterization of the continuous functions on $[-1,1]$ which result from the restriction of an analytic function on the interior of the ellipse
$$E_{\rho}=\{(\frac{1}{2}(\rho+\frac{1}{\rho})\cos\theta, \frac{1}{2}(\rho-\frac{1}{\rho})\sin\theta):\theta\in [0,2\pi]\},$$
as those functions $f\in C[-1,1]$ such that $\lim\sup_{n\to\infty}E(f,\Pi_n)^{\frac{1}{n}}=\frac{1}{\rho}$.

In this short note we prove that, if $(C[a,b],\{A_n\})$ is an approximation scheme and $(A_n)$ satisfies de La Vallée-Poussin Theorem, there are instances of continuous functions on $[a,b]$,
real analytic on $(a,b]$, which are ``poorly approximable'' by the elements of the approximation scheme $(A_n)$. This illustrates the thesis
that the smoothness conditions guaranteeing that a function is ``well approximable''
must be, at least in these cases, ``global''. The failure of smoothness at endpoints may result in an
arbitrarily slow rate of approximation.  A result of this kind,  which is highly nonconstructive, based on different arguments, and applicable to different approximation schemes, was  recently proved by Almira and Oikhberg \cite{almira_oikhberg}.

\section{The main result}
Let us start by recalling what an approximation scheme is.
\begin{definition}[See, for example, \cite{almira}, \cite{almiraluther2}]
Let $(X,\|\cdot\|)$ be a quasi-Banach space, and let
$\{0\}=A_0\subset A_1\subset\cdots \subset A_n\subset\cdots \subset X$
be an infinite chain of subsets of $X$, where all inclusions are
strict. We say that $(X,\{A_n\})$ is an {\it approximation scheme} (or that $(A_n)$ is an approximation scheme in $X$) if:
\begin{itemize}
\item[$(i)$] There exists a map $K:\mathbb{N}\to\mathbb{N}$ such that $K(n)\geq n$ and $A_n+A_n\subseteq A_{K(n)}$ for all $n\in\mathbb{N}$.

\item[$(ii)$] $\lambda A_n\subset A_n$ for all $n\in\mathbb{N}$ and all scalars $\lambda$.

\item[$(iii)$] $\bigcup_{n\in\mathbb{N}}A_n$ is a dense subset of $X$
\end{itemize}
A particular example is a {\it linear approximation scheme}, arising when the
sets $A_n$ are linear subspaces of $X$. In this setting, we can take $K(n) = n$.
An approximation scheme is called {\it non-trivial} if $X \neq \cup_n \overline{A_n}$.
\end{definition}

\begin{definition} The vector space $A \subset C[a,b]$ is named Haar on $[a,b]$ if $\dim A=n$ and the only element from $A$ which has more than $n-1$ zeroes is the null function.
\end{definition}

\begin{definition} Let $(C[a,b],\{A_n\})$ be an approximation scheme. We say that $(A_n)$ satisfies de La Vallée-Poussin Theorem if there exists a function $\phi_{\{A_n\}}:\mathbb{N}\to\mathbb{N}$, $\phi_{\{A_n\}}(n)=m_n$, such that, if $f\in C[a,b]$ and $a\leq t_{n,0}<t_{n,1}<\cdots <t_{n,m_n}\leq b$ satisfy  $\mathbf{sign}(f(t_{n,k})f(t_{n,k+1}))=-1$ for all $k$ and $\min_{0\leq k<m_{n}}|f(t_{n,k})|>\varepsilon_n$, then $E(f,A_n)\geq \varepsilon_n$.  \end{definition}

Classical de La Vallée-Poussin Theorem \cite{poussin} (see also \cite[page 74, Theorem 5.2]{devore}) appears when $(C[a,b],\{A_n\})$ is a linear approximation scheme, $m_n=\dim A_n$,  and each $A_n$ is Haar. More examples of approximation schemes satisfying de La Vallée-Poussin Theorem appear in connection with rational approximation. Indeed, $\{n_k\},\{m_k\}\to +\infty$ are a pair of increasing sequences of natural numbers, and $R_{n,m}$ denotes the set of rational functions $r(t)=p(t)/q(t)$ with poles outside $[a,b]$ and $\deg(p)\leq n$, $\deg(q)\leq m$, the approximation scheme $(C[a,b],\{R_{n_k,m_k}\}_{k=0}^\infty)$ satisfies de La Vallée-Poussin Theorem \cite[Theorem 98]{meinardus}. Finally, it is an easy exercise to check that if the approximation scheme $(C[a,b],\{A_n\})$ satisfies de La Vallée-Poussin Theorem with $\phi_{\{A_n\}}(n)=m_n$ and we denote by $\Sigma_n(A_n)$ the set of splines with $n$ free knots constructed with the elements of $A_n$ then $(C[a,b],\{\Sigma_n(A_n)\})$  satisfies de La Vallée-Poussin Theorem with $\phi_{\{\Sigma_n(A_n)\}}(n)=(n+1)^2m_n$. Recall that  $f\in \Sigma_n(A_n)$ if and only if there exists a partition $a=t_0<t_1<\cdots<t_n<b=t_{n+1}$ such that  $f(t)=\sum_{k=0}^{n-1} a_k(t)\chi_{[t_k,t_{k+1})}(t)+a_k(t)\chi_{[t_n,t_{n+1}]}(t)$ for some $a_0,a_1,\cdots a_n\in A_n$. Here, $\chi_{I}$ denotes the characteristic function associated to the interval $I$.

\begin{theorem}\label{mainresult} Let us assume that $0<\alpha<\beta$ and let $[a,b]=[0,1]$ or $[a,b]=[\alpha,\beta]$. Let $(C[a,b],\{A_n\})$ be an approximation scheme such that  $(A_n)$ satisfies de La Vallée-Poussin Theorem with $\phi_{\{A_n\}}(n)=m_n$. Let $\{\varepsilon_n\}_{n=0}^{\infty}\in c_0$ be a non-increasing sequence of positive numbers converging to zero. Then there exists $f\in C[a,b]$ such that $f$ is real analytic on $(a,b]$ and $E(f,A_n)\geq \varepsilon_n$ for all $n=0,1,2,\cdots$.
\end{theorem}
\noindent \textbf{Proof.} Let us first consider the case $[a,b]=[0,1]$.  The result is trivial if $\{\varepsilon_n\}$ is stationary, so that we can assume $\varepsilon_n>0$ for all $n\in\mathbb{N}$. Let $p:\mathbb{R}\to\mathbb{R}$ be the continuous polygonal line with vertices $\{(n,p(n))\}_{n\in\mathbb{Z}}$ defined by
  \begin{itemize}
  \item[$(i)$] $p(t)$ is an even function, $p(0)=p(1)=\cdots=p(m_0-1)=3\varepsilon_0$.
  \item[$(ii)$] $p(m_0+m_1+\cdots+m_n+k)=3\varepsilon_{n+1}$ for all $n\in\mathbb{N}$ and $k=0,1,\cdots,m_{n+1}-1$.
  \end{itemize}
It follows from Carleman's Theorem \cite{carleman} (see also \cite[Chapter 1, Theorem 4.3]{makovoz}) that there exists an entire function $e(t)$ such that $e(\mathbb{R})\subset \mathbb{R}$ and $\sup_{t\in\mathbb{R}}|p(t)-e(t)|\leq \frac{p(t)}{3}$. Hence $e([0,m_0-1]) \subseteq (\varepsilon_0,5\varepsilon_0)$ and, for $n\in\mathbb{N}$, $e([m_0+\cdots+m_n,m_0+\cdots+m_{n+1}-1))\subseteq (\varepsilon_{n+1},5\varepsilon_{n+1})$. Let us set $f(0)=0$ and, for $t\in (0,1]$, $f(t)=e(\frac{1}{t})\cos(\frac{2\pi}{t})$. This is our function. Obviously, $f$ is analytic on $(0,1]$ and continuous on $[0,1]$. To prove that $E(f,A_n)\geq \varepsilon_n$ it is enough to take into account that $A_n$ satisfies de La Valleé-Poussin Theorem with $\phi_{\{A_n\}}(n)=m_n$ and, by construction, there exists an ordered set of points
$\frac{1}{m_0+m_1+\cdots+m_{n-1}}\leq t_{n,0}<t_{n,1}<\cdots <t_{n,m_n}\leq \frac{1}{m_0+m_1+\cdots+m_{n}-1}$ such that $\mathbf{sign}(f(t_{n,k})f(t_{n,k+1}))=-1$ for all $k$ and $\min_{0\leq k<m_{n}}|f(t_{n,k})|>\varepsilon_n$. This ends the proof for $[a,b]=[0,1]$. If $0<a<b$, the arguments above show that $E(g,A_n)\geq \varepsilon_n$ for all $n$, where $g(t)=f(\frac{t-a}{b-a})$. {\hfill $\Box$}

\begin{remark} The same kind of phenomenon also holds in different situations for higher dimensions. To prove this claim, let $(C([-1,1]^s),\{A_n\}_{n=0}^{\infty})$ be an approximation scheme such that, for each $k\in\{1,\cdots,s\}$ and $a=(a_1,\cdots,a_s)\in [-1,1]^s$, the sets $B_{k,n}^a=\{g(x)=f(a_1,\cdots,a_{k-1},x,a_{k+1},\cdots,a_s):f\in A_n\}$ form an approximation scheme in $C[0,1]$ and $(C[0,1],\{B_{k,n}^a\}_{n=0}^\infty)$ satisfies de La Vallée-Poussin Theorem with $\phi_{\{B_{k,n}^a\}}(n)=m_{n}$. These approximations schemes are easy to find. For example, for $s=2$ we can take $$A_n=\Pi_{n}[x,y]:=\{P(x,y)=\sum_{0\leq i,j\leq i+j\leq n}a_{ij}x^iy^j:a_{ij}\in\mathbb{R} \text{ for all }i,j\}$$ or $$A_n=\{P(x,y)/Q(x,y):P,Q\in \Pi_n[x,y] \text{ and } Q(a,b)\neq 0 \text{ for all } (a,b)\in [-1,1]^2\}.$$ Then, given a non-increasing sequence $\{\varepsilon_n\}\in c_0$ we choose $f\in C[0,1]$ as in Theorem \ref{mainresult} and define $g(x_1,\cdots,x_s)=f(x_1^2+x_2^2+\cdots+x_s^2)$. This function is continuous on $[-1,1]^s$ and real analytic on $\Omega=[-1,1]^s\setminus\{\mathbf{0}\}$, where $\mathbf{0}=(0,0,\cdots,0)\in \mathbb{R}^s$ (see \cite[Proposition 2.2.8]{krantz}). Moreover, $h(x)=f(x^2)$ also satisfies the conclusion of Theorem \ref{mainresult} since the change of $x$ by the new variable $x^2$ does not eliminates the oscillations of $f$. Hence, for each $n\in\mathbb{N}$ and $k\in\{1,\cdots,s\}$, we have that $E(g,A_n)_{C([-1,1]^s)}\geq E(g(0,\cdots,0,x_k,0,\cdots,0),B_{k,n}^{\mathbf{0}})_{C[0,1]}=E(h,B_{k,n}^{\mathbf{0}})_{C[0,1]}\geq \varepsilon_n$. \end{remark}

%%%  ==============================================================

\bigskip

\footnotesize{J. M. Almira

Departamento de Matem\'{a}ticas. Universidad de Ja\'{e}n.

E.P.S. Linares,  C/Alfonso X el Sabio, 28

23700 Linares (Ja\'{e}n) Spain

email: jmalmira@ujaen.es}


\begin{thebibliography}{13}
\bibitem{almira}   \textbf{J. M. Almira, U. Luther, } Inverse closedness of approximation algebras. J. Math. Anal. Appl. \textbf{314} (2006), no. 1, 30--44.

\bibitem{almiraluther2} \textbf{J. M. Almira, U. Luther, } Generalized approximation spaces and applications, Math. Nachr. \textbf{263-264} (2004) 3--35.

 \bibitem{almira_oikhberg} \textbf{J. M. Almira, T. Oikhberg, } Shapiro's Theorem for subspaces,  J. Math. Anal. Appl. (2011) , doi:10.1016/j.jmaa.2011.09.054
%\bibitem{carleson} \textbf{L. Carleson, } An interpolation problem for bounded analytic functions, Amer. J. of Math. \textbf{80} (1958) 921-930.

\bibitem{bernsteininverso}  \textbf{S. N. Bernstein, } Sur l'ordre de la meilleure approximation des fonctions continues par les polynômes de degré donné, Mem. Cl. Sci. Acad. Roy. Belg. \textbf{4} (1912), 1-103.

\bibitem{carleman} \textbf{T. Carleman}   Sur un théorème de Weierstrass,  Arkiv. Mat. Astron. Fys., \textbf{20} (1927) 1-5.

\bibitem{poussin} \textbf{C. de la Vallée Poussin, } \textit{Lessons sur l'Approximation des Fonctions d'une Variable Réelle, } Paris, Gauthier-Villars, 1919.

\bibitem{devore} \textbf{R. A. DeVore, G. G. Lorentz, } \textit{Constructive approximation, } Springer, 1993.

%\bibitem{hoffman} \textbf{K. Hoffman, } \textit{Banach spaces of analytic functions, } Dover, 1988.

%\bibitem{hayman} \textbf{W. K. Hayman, } Interpolation by bounded functions, Annal. del Institut Fourier
%\textbf{8} (1958), 277-290.

\bibitem{jackson} \textbf{D. Jackson,  } Über die Genauigkeit der Annäherung stetiger Funktionen durch ganze rationale Funktionen gegebenen Grades und trigonometrische Summen gegebener Ordnung (On the precision of the approximation of continuous functions by polynomials of given degree and by trigonometric sums of given order). Preisschrift und Dissertation. Univ. Göttingen, June 14, 1911 (at GDZ). This is Dunham Jackson's doctoral thesis.

\bibitem{krantz}  \textbf{S. G. Krantz, H. R. Parks, } \textit{A primer of Real Analytic Functions. } Birkh\"{a}user, 2002.

\bibitem{makovoz} \textbf{G. G. Lorentz, M. von Golitschek, Y. Makovoz, } \textit{Constructive approximation. Advanced problems. } Springer, 1996.

\bibitem{meinardus} \textbf{G.~Meinardus, } \textit{Approximation of functions: theory and numerical methods, } Springer-Verlag, 1967.


\bibitem{pinkus} \textbf{A. Pinkus, } Negative Theorems in Approximation Theory, American Math. Monthly \textbf{110} (2003), 900-911.

\end{thebibliography}
\end{document}